\documentclass{amsart}

\newtheorem{theorem}{Theorem}

\theoremstyle{definition}
\newtheorem{Definition}[theorem]{Definition}

\usepackage{amsfonts}

\newcommand{\fF}{{\mathfrak F}} 

\newcommand{\id}{\triangleleft}

\newcommand{\sN}{\mathcal{N}} 
\newcommand{\cent}{\mathcal{C}}

\DeclareMathOperator{\core}{Core}

\begin{document}

\title{Lie algebra $\fF$-normalisers are intravariant}

\author{Donald W. Barnes}
\address{1 Little Wonga Rd. Cremorne NSW 2090 Australia}
\email{donwb@iprimus.com.au}
\thanks{This work was done while the author was an Honorary Associate of the School
of Mathematics and Statistics, University of Sydney}

\subjclass[2000]{Primary 17B30}
\keywords{Soluble Lie algebras, formations, intravariant}
\date{}

\begin{abstract}  Let $\fF$ be a saturated formation of soluble Lie algebras.  Let $L$ be a soluble
Lie algebra and let $U$ be an $\fF$-normaliser of $L$.  Then $U$ is intravariant in $L$.
\end{abstract}

\maketitle

All Lie algebras considered in this paper are finite-dimensional over the field $F$ and are soluble.  The
theory of saturated formations and of $\fF$-projectors was developed in Barnes and Gastineau-Hills
\cite{BGH}.  The theory of $\fF$-normalisers was developed in Stitzinger \cite{Stit}.   The following
definition of intravariance for subalgebras of Lie algebras was given in Barnes \cite{Frat}.

\begin{Definition}  The subalgebra $U$ of the Lie algebra $L$ is said to be
{\em intravariant} in $L$ if every derivation of $L$ is expressible as the sum of an inner
derivation and a derivation which stabilises $U$.  \end{Definition}

By \cite[Lemma 1.2]{Frat}, the intravariant subalgebras of $L$ are precisely those subalgebras $U$ with
the property that, if $L$ is an ideal of $L^*$, then $L^* = L + \sN_{L^*}(U)$.  It was proved in 
\cite[Theorem 2.2]{Frat} that, if $\fF$ is a saturated formation of soluble Lie algebras, then the
$\fF$-projectors of a soluble Lie algebra $L$ are intravariant in $L$.  It seems reasonable to
conjecture that the $\fF$-normalisers also are intravariant.  The following definitions are taken from
Stitzinger \cite{Stit}.

\begin{Definition}  A maximal subalgebra $M$ of $L$ is called $\fF$-normal if it
complements an $\fF$-central chief factor of $L$.  Equivalently, $M$ is
$\fF$-normal if for some chief factor $A/B$ of $L$, $M+A = L$ and the split
extension of $A/B$ by $M/\cent_M(A/B)$ is in $\fF$.  This is equivalent to
$L/\core(M) \in \fF$.  $M$ is called $\fF$-abnormal if it is not $\fF$-normal, that
is, if $L/\core(M) \notin \fF$. \end{Definition}

\begin{Definition} A maximal subalgebra $M$ of $L$ is called $\fF$-critical if $M$
is $\fF$-abnormal and $M + N(L) = L$ where $N(L)$ is the nil radical of $L$. \end{Definition}

\begin{Definition} A subalgebra $V \le L$ is called an $\fF$-normaliser of $L$ if
$V \in \fF$ and there exists a chain of subalgebras
$$L = M_0 > M_1 > \ldots > M_n = V$$
with each $M_i$ $\fF$-critical in $M_{i-1}$. \end{Definition}

We shall need the following result of Stitzinger \cite{Stit}.

\begin{theorem} Let $U$ be an $\fF$-normaliser of $L$.  Then $U$ covers every
$\fF$-central chief factor and avoids every $\fF$-eccentric chief factor of $L$.\end{theorem}

We now prove the asserted result.

\begin{theorem}\label{th-main} Let $\fF$ be a saturated formation.  Let $L$ be a soluble Lie
algebra and let $U$ be an $\fF$-normaliser of $L$.  Then $U$ is intravariant in $L$.
\end{theorem}

\begin{proof}   Let $d$ be a derivation of $L$ and let $D = \langle d, L\rangle$ be the split extension of
$L$ by $d$.  We have to prove that $\sN_D(U) + L = D$.  We use induction over $\dim{L}$.  Let $A
\subseteq L$ be a minimal ideal of $D$.  Then $U+A/A$ is an $\fF$-normaliser of $L/A$, so by induction,
$\sN_D(U+A) + L = D$.  Let $N = \sN_D(U+A)$.  Since $A$ is an irreducible $D$-module and $L$ is an ideal
of $D$, all $L$-composition factors of $A$ are isomorphic, so either all are $\fF$-central and $U
\supset A$ or all are $\fF$-eccentric and $U \cap A = 0$.  If $U \supset A$, the $\sN_D(U) = N$ and the
result holds.  If $U \cap A = 0$, then $U$ is an $\fF$-projector of $U+A$ and is intravariant in $U+A$. 
As $U+A \id N$, we have $\sN_N(U) + (U+A) = N$ and $\sN_D(U) + L \subseteq  \sN_N(U) + (U+A) + L =
N+L = D$.
\end{proof}

\bibliographystyle{amsplain}

\begin{thebibliography}{10}

\bibitem{Frat} D. W. Barnes, \textit{The Frattini argument for Lie algebras}, J.
Algebra  \textbf{27} (1973), 486--490.

\bibitem{BGH} D. W. Barnes  and H. M. Gastineau-Hills, \textit{On the theory of
soluble Lie algebras}, Math. Zeitschr. \textbf{106} (1968), 343--354.

\bibitem{Stit} E. L. Stitzinger, \textit{Covering-avoidance for saturated
formations of solvable Lie algebras}, Math. Zeitschr. \textbf{124} (1972), 237--249.


\end{thebibliography}

\end{document}